\documentstyle[12pt]{article}

\newcommand{\C}{{\bf C}}
\newcommand{\R}{{\bf R}}

\newcommand{\Z}{{\bf Z}}

\begin{document}

\begin{center}
{\large \bf Dimers and Dominoes}
\end{center}

\begin{center}
{\sc September 24, 1992} \\
(revised November 4, 1997)
\end{center}

\vspace{5ex}
\begin{center}
James Propp \\
Massachusetts Institute of Technology and MSRI \\[2ex]
\end{center}

\vspace{0.2in}

\noindent
{\bf 1. Counting configurations on a rectangle.}

\vspace{0.1in}

Let us count the number of domino tilings of an $m$-by-$n$ checkerboard,
or equivalently, the number of matchings
of an $m$-by-$n$ square grid of dots, $G$,
where by a {\bf matching} of a graph we mean a collection of edges
such that each vertex of the graph belongs
to exactly one edge in the collection.
(In physics, one regards each vertex as an atom,
and each edge in a matching as representing
a diatomic molecule, or dimer;
hence a matching of a graph
is also known as a dimer cover.)
In order for the checkerboard to have a domino tiling,
it is necessary that the number of cells, $mn$, be even;
assume for definiteness that $n$ is even.

Let $N(m,n)$ denote the number of tilings.
We know that $N(m,n)$ must be at least $(F_m)^{n/2}$,
where $F_m$ denotes the $m$th Fibonacci number,
because (as is easily shown by induction)
there are $F_m$ ways of tiling a 2-by-$m$ rectangle
with dominoes.
That is, $N(m,n)$ is asymptotically
at least as large as $\sqrt{\phi}^{mn}$,
where $\phi=(1+\sqrt{5})/2=1.62\cdots$
and $\sqrt{\phi}=1.27\cdots$.
On the other hand, we can encode any tiling of the region
with at most $mn/2$ bits, one tile at a time,
by iterating the following rule:
Find the upper-left-most cell
that does not belong to an already-encoded domino,
and record 0 or 1
according to whether it shares a domino
with the cell to its right or the cell below it.
For the price of one bit we have encoded the position
of a domino covering two cells,
and so it takes only $mn/2$ bits
to encode the full tiling.
Thus $N(m,n)$ is asymptotically
no larger than $\sqrt{2}^{mn}$.

We will find it more convenient to work with
matchings of $G$ than tilings of the checkerboard,
but we want to retain the notion of an alternating coloring.
We color the vertices of $G$ black and white
such that each black vertex has only white neighbors and vice versa.

Let $A$ be the adjacency matrix of the graph $G$,
whose rows and columns are indexed by the vertices of $G$
and whose $j,k$th entry equals 1 if the $j$th vertex of $G$
is adjacent to the $k$th vertex of $G$
and equals 0 otherwise.
Let $B$ be the ``bipartite adjacency matrix'' of $G$,
whose rows and columns are indexed by
the black vertices and white vertices of $G$, respectively,
and whose $j,k$th entry equals 1 if the $j$th black vertex of $G$
is adjacent to the $k$th white vertex of $G$
and equals 0 otherwise.
If we order the vertices of $G$
so that all the black vertices precede all the white vertices,
then $A$ has the block-structure
\[ 
A = \left( \begin{array}{cc} 0 & B \\ B^{\mbox t} & 0 \end{array} \right)
\]
(where the superscript ``t'' denotes the transpose).
It is clear that the permanent of $A$
is equal to the square of the permanent of $B$,
and that the permanent of $B$
is equal to the number of matchings of $G$.

Our approach is based on that of Kasteleyn
(P.W. Kasteleyn, {\it The statistics of dimers on a lattice, I:
The number of dimer arrangements on a quadratic lattice,}
Physica {\bf 27} (1961), 1209-1225).
Let $\tilde{A}$ and $\tilde{B}$
be the matrices obtained from $A$ and $B$, respectively,
by replacing each entry $+1$ that corresponds to a vertical bond in $G$ by $+i$.
I claim that $\det \tilde{B} = \pm N(m,n)$,
where $N(m,n)$ is the number of matchings of the graph.
To prove this, it suffices to show that
all the contributions to the determinant
equal $\pm 1$, and that all have the same sign.
If all bonds are horizontal,
the contribution is clearly $\pm 1$,
since no $i$'s are involved.
We know (see, for instance, W. Thurston, {\it Conway's tiling groups,}
Amer. Math. Month. 97 (1990), 757-773)
that every matching can be obtained
from the all-horizontal matching by means of
elementary moves of the form
$$\cdots (a \leftrightarrow c) (b \leftrightarrow d) \cdots
\ \Leftrightarrow \ 
\cdots (a \leftrightarrow d) (b \leftrightarrow c) \cdots \ .$$
In the original matrix $B$, this corresponds to a transposition
\[
\left( \begin{array}{cc} 1 & 0 \\ 0 & 1 \end{array} \right)
\leftrightarrow
\left( \begin{array}{cc} 0 & 1 \\ 1 & 0 \end{array} \right)
\]
(where the two rows are indexed $a$ and $b$,
the two columns are indexed $c$ and $d$,
and all intervening rows and columns are omitted),
which yields a transversal of the matrix whose 
parity is opposite that of the original transversal.
But in the matrix $\tilde{B}$,
this change of parity is compensated for by the fact that
in one of the two transversals,
the associated entries are $+i$'s instead of $+1$'s.
Note that
\[ 
\tilde{A} = \left( \begin{array}{cc} 0 & \tilde{B} \\ 
\tilde{B}^{\mbox t} & 0 \end{array} \right) , \] so $\det\tilde{A} = \pm (\det\tilde{B})^2$.
To determine the sign, note that
if we combine any transversal of $\tilde{B}$
with its own transpose,
we get a transversal of $\tilde{A}$
that is the product of an even number of transpositions.
This transversal contributes $+1$ to $\det\tilde{A}$,
so $\det\tilde{A} > 0$,
and $\det\tilde{A} = N(m,n)^2$.

Here is another way of verifying
that the terms in the determinant
add coherently,
i.e. with no cancellations.
Label the white nodes of $G$ from 1 to $mn/2$,
and do the same for the black nodes,
so that every matching of $G$
corresponds to a permutation of $\{1,...,mn/2\}$
and can be assigned a parity,
which coincides with the parity of the associated transversal of
the matrix $B$.
Put $+1$'s on the horizontal edges of the graph $G$
and $+i$'s on the vertical edges,
and call these the {\bf weights} of the edges.
Say more generally that the weight of a set of edges
is the product of the weights of the edges belonging to the set.
Every matching $\mu$ of $G$ 
corresponds to a non-vanishing transversal of $\tilde{B}$,
and its contribution to $\det\tilde{B}$
is equal to its parity $\sigma(\mu)$ 
times its weight $w(\mu)$,
where $w(\mu)$ must equal $\pm 1$
since every matching contains an even number of vertical edges.
We need to show that for any two matchings $\mu_1$ and $\mu_2$ of $G$,
$\sigma(\mu_1) w(\mu_1) = \sigma(\mu_2) w(\mu_2)$,
or equivalently,
$\sigma(\mu_1) / \sigma(\mu_2) = w(\mu_2) / w(\mu_1)$.

Let $\mu_1 + \mu_2$ be the multigraph obtained by
combining the edges of $\mu_1$ and the edges of $\mu_2$.
Thus, every edge that is common to $\mu_1$ and $\mu_2$
yields a double edge (which we regard as a 2-cycle),
while the rest of $\mu_1 + \mu_2$
splits up into even-length cycles of length 4 or more.
Suppose we get $r$ cycles in this way;
then $\mu_2$ can be obtained from $\mu_1$
in $r$ stages,
by rotating edges around cycles.
represent their lengths by $2l_1, 2l_2, ..., 2l_r$,
where $2l_1+2l_2+...+2l_r=mn$.
The permutation of $\{1,...,mn/2\}$ associated with $\mu_2$
is equal to
the permutation associated with $\mu_1$
times the product of an $l_1$-cycle, an $l_2$-cycle, ... , an $l_r$-cycle.
Thus the parity of $\mu_2$
equals the parity of $\mu_1$
times $(-1)^{(l_1+1)+(l_2+1)+...+(l_r+1)}$.
(Note: the length of the graph-theoretic cycle in $G'$
is twice the length of the associated permutational cycle on $\{1,...,mn/2\}$.)
That is, a cycle in $\mu_1 + \mu_2$ of length $2l$ 
contributes a factor of $(-1)^{l+1}$
to the relative parity of $\mu_2$ with respect to $\mu_1$.

What about the relative weight of $\mu_2$ with respect to $\mu_1$?
Given an oriented closed path $\Gamma$ in the square grid,
define the {\it disparity} $d(\Gamma)$
as the number of vertical edges in $\Gamma$
that point from a black to a white vertex
minus the number of vertical edges in $\Gamma$
that point from a white to a black vertex.
Then it can be seen that
the contribution to the relative weight $w(\mu_2) / w(\mu_1)$
made by a particular cycle in $\mu_1 + \mu_2$
is equal to $i^{d(\Gamma)}$,
where $\Gamma$ is an (arbitrary) orientation of that cycle.
Since $\Gamma$ is non-self-intersecting
it is easy to show that
$d(\Gamma)$ is equal to the sum of $d(\Gamma_c)$
as $c$ runs over the square cells enclosed by $\Gamma$,
where $\Gamma_c$ is the closed path of length 4
that encircles the $j$th square cell enclosed by $\Gamma$
and has the same orientation as $\Gamma$ itself;
for if one examines the sum,
one sees that cancellation takes place
on all edges interior to $\Gamma$
that belong to two $\Gamma_c$'s,
leaving only edges on $\Gamma$ itself uncancelled.
But $d(\Gamma_c)$ is $\pm 2$,
according to whether $c$ is a black or white cell
under a checkerboard coloring.
Hence $d(\Gamma)$ is equal to twice the difference
between the number of black cells and the number of white cells
enclosed by $\Gamma$,
which is congruent, modulo 4,
to twice the sum
of the number of black cells and the number of white cells
enclosed by $\Gamma$.
Hence a cycle in $\mu_1 + \mu_2$
enclosing area $A$
contributes a factor of $(-1)^{A}$
to the relative weight of $\mu_2$ with respect to $\mu_1$.

To finish our analysis,
recall Pick's theorem for polygons 
whose vertices belong to a square grid:
$A = I + \frac{1}{2}B - 1$,
where $A$ denotes the area enclosed by a polygon,
$I$ denotes the number of interior grid-points,
and $B$ denotes the number of grid-points on the boundary.
In the case of our cycle $\Gamma$,
the number of interior vertices $I$
must be even
(since both $\mu_1$ and $\mu_2$
give matchings on this set of vertices),
and $\frac{1}{2}B$ is just $l$,
the length of the cycle.
Hence $A$ is congruent to $l-1$ modulo 2,
so that the parity-factor $(-1)^{l+1}$
and the weight-factor $(-1)^A = (-1)^{l-1}$
exactly cancel.

Thus we see that all the matchings
contribute coherently to the determinant of $B$.

Example: $m=3$, $n=2$.  Index the entries in a rectangle
with 3 rows and 2 columns as shown:
\[ \begin{array}{cc} 1 & 4 \\ 5 & 2 \\ 3 & 6 \end{array} \]
Then we get
\[ \det\tilde{A} = \left| \begin{array}{cccccc}
	0 & 0 & 0 & 1 & i & 0 \\
	0 & 0 & 0 & i & 1 & i \\
	0 & 0 & 0 & 0 & i & 1 \\
	1 & i & 0 & 0 & 0 & 0 \\
	i & 1 & i & 0 & 0 & 0 \\
	0 & i & 1 & 0 & 0 & 0 \end{array} \right| \ = \ 9 = 3^2, \]
corresponding to the fact that the 3-by-2 rectangle
can be tiled by dominoes in exactly 3 ways.

We will find a general formula for $\det\tilde{A}$
by determining the full spectrum of $\tilde{A}$.
Here I use an approach suggested to me by Noam Elkies
in private correspondence.
Let $V$ be the space of functions $f: \Z^2 \rightarrow \C$
such that $f(x,-y)=f(-x,y)=-f(x,y)$
and $f(x+2(m+1),y)=f(x,y+2(n+1))=f(x,y)$ for all $x,y \in \Z$.
Note that for such a function $f$,
$f(x,y)$ vanishes if either
$x$ is a multiple of $m+1$
or $y$ is a multiple of $n+1$.
Let $L: V \rightarrow V$ be the modified local summation operator
$(Lf)(x,y) = f(x-1,y) + f(x+1,y) + i f(x,y-1) + i f(x,y+1)$.
It is easy to check that $L$ does indeed send $V$ to $V$.
$V$ is isomorphic to $\C^{mn}$
(restrict $f:\Z^2 \rightarrow \C$ to $\{1,...,m\} \times \{1,...,n\}$),
and the $mn$-by-$mn$ matrix that intertwines the action of $L$
under this automorphism is $\tilde{A}$.
Hence $\det\tilde{A}$ is the product of the eigenvalues of $L$.

A basis for $V$ is given by the functions
\[ f_{j,k} (x,y) = \sin \frac{\pi j x}{m+1} \sin \frac{\pi k y}{n+1} \]
($1 \leq i \leq m$, 
$1 \leq j \leq n$). 
It is easy to check that these are eigenfunctions of $L$,
using the identity 
$\sin (\alpha-\beta) + \sin (\alpha+\beta) = (\sin \alpha) (2 \cos \beta)$:
\begin{eqnarray*}
(Lf_{j,k})(x,y) & = &
	\left(\sin \frac{\pi j (x-1)}{m+1} + 
		\sin \frac{\pi j (x+1)}{m+1}\right)
	\sin \frac{\pi k y}{n+1} \\
& & + \ i \sin \frac{\pi j x}{m+1}
	\left(\sin \frac{\pi k (y-1)}{n+1} + 
		\sin \frac{\pi k (y+1)}{n+1}\right) \\
& = &	\left(\sin \frac{\pi j x}{m+1}\right) 
	\left(2 \cos \frac{\pi j}{m+1}\right)
	\left(\sin \frac{\pi k y}{n+1}\right) \\
& & + \ i \left(\sin \frac{\pi j x}{m+1}\right)
	\left(\sin \frac{\pi k y}{n+1}\right) 
	\left(2 \cos \frac{\pi k}{n+1}\right) \\
& = & \left(2 \cos \frac{\pi j}{m+1} + 2i \cos \frac{\pi k}{n+1}\right)
	(Lf_{j,k}) (x,y) \ .
\end{eqnarray*}
So
\[ \det\tilde{A} = \prod_{j=1}^m \prod_{k=1}^n
\left(2 \cos \frac{\pi j}{m+1} + 2i \cos \frac{\pi k}{n+1}\right) . \]
Since $n$ is even, we can combine the $k$ and $n+1-k$ factors to get
\[ \prod_{j=1}^{m} \prod_{k=1}^{n/2}
\left(4 \cos^2 \frac{\pi j}{m+1} + 4 \cos^2 \frac{\pi k}{n+1}\right) . \]
Since the $j$ and $m+1-j$ factors are equal,
and since we want the positive square root of $\det\tilde{A}$ anyway,
we get
\[ N(m,n) = \prod_{j=1}^{m/2} \prod_{k=1}^{n/2}
\left(4 \cos^2 \frac{\pi j}{m+1} + 4 \cos^2 \frac{\pi k}{n+1}\right) \]
when $m$ is even.

Thus, for instance, the number of domino tilings of an 8-by-8 checkerboard
is $12988816 \approx 1.29^{8^2}$.
(Note that 12988816 is exactly equal to $3604^2$.
In fact, $N(n,n)$ is a perfect square when $n$ is congruent to 0 mod 4
and twice a perfect square when $n$ is congruent to 2 mod 4.
This is a special case of a result of William Jockusch's;
he showed that if $G$ is a bipartite graph
with a 4-fold rotational symmetry that swaps the two color classes,
then the number of matchings of $G$
is a square if the the number of vertices of $G$ is a multiple of 8
and twice a square otherwise.)

Another way to calculate $\det\tilde{A}$
is to note that it is of the form
$C_m \otimes I_n + I_m \otimes iC_n$
where 
$\otimes$ denotes the Kronecker product,
$I_m$ is the $m$-by-$m$ identity matrix,
$I_n$ is the $n$-by-$n$ identity matrix,
$C_m$ is the $m$-by-$m$ matrix
\[ \left( \begin{array}{ccccccc}
0 & 1 & 0 & 0 & \dots & 0 & 0 \\
1 & 0 & 1 & 0 & \dots & 0 & 0 \\
0 & 1 & 0 & 1 & \dots & 0 & 0 \\
\vdots & \vdots & \vdots & \vdots & & \vdots & \vdots \\
0 & 0 & 0 & 0 & \dots & 0 & 1 \\
0 & 0 & 0 & 0 & \dots & 1 & 0 \end{array} \right) \ , \]
and $C_n$ is the $n$-by-$n$ matrix of the same form.
The eigenvalues of $\det\tilde{A}$
are therefore precisely the values $\mu + i \nu$
where $\mu$,$\nu$ are eigenvalues of $C_m$ and $C_n$, respectively.
This yields the same answer as before.

Kasteleyn's approach is slightly different:
instead of using a determinant, he uses a Pfaffian.
However, since our graph is bipartite,
it's not hard to show that the two methods
are algebraically equivalent.

In general, for $n$ even we have
\[ \log N(n,n) = \frac{1}{2} \sum_{j=1}^n \sum_{k=1}^n
	\ \log \, \left(2 \cos \frac{\pi j}{n+1} 
		+ 2i \cos \frac{\pi k}{n+1}\right) \]
so that as $n$ gets large
\begin{eqnarray*}
\frac{1}{n^2} \log N(n,n) & \rightarrow &
\frac{1}{2} \int_0^1 \int_0^1 \log \ (2 \cos \pi s + 2i \cos \pi t) \ ds \: dt
\end{eqnarray*}
which is one-half of the average value of
$\log \: (\alpha + \alpha^{-1} + i \beta + i \beta^{-1})$
as $\alpha$ and $\beta$ range independently 
over the unit circle with uniform density.

It can be shown that the double integral evaluates to $G/\pi$, where $G$ is
Catalan's constant $1/1 - 1/9 + 1/25 - 1/49 + ...$.
Thus, the number of domino tilings of the $n$-by-$n$ board (with $n$ even)
is roughly $1.34^{n^2}$.

\vspace{0.2in}

{\bf 2. Counting configurations on a torus.}

\vspace{0.1in}

Let us now count the number of domino tilings of an $m$-by-$n$ torus,
or equivalently, the number of matchings
of an $m$-by-$n$ toroidal grid of dots, $G'$.
$G'$ can be obtained from $G$
by adding $m$ extra horizontal edges
and $n$ extra vertical edges;
we call these edges of $G'$ {\bf special}.
For simplicity, we will require that both $m$ and $n$ be even.

It turns out, for reasons that involve the non-planarity of the graph $G'$,
that there is no way of replacing the entries of the adjacency matrix
by roots of unity to obtain a new matrix 
whose determinant equals the permanent of the original adjacency matrix.
However, we will see that
the number of matchings of $G'$,
which we denote by $N'(m,n)$,
can be written as a linear combination of four determinants.

Let $A'$ and $B'$ be
the adjacency matrix
and bipartite adjacency matrix of $G'$, respectively,
so that 
\[ 
A' = \left( \begin{array}{cc} 0 & B' \\ B'^{\mbox t} & 0 \end{array} \right) .
\]
Let $\tilde{A'}$ and $\tilde{B'}$
be the matrices obtained from $A$ and $B$, respectively, by replacing 
each entry $+1$ that corresponds to a vertical bond in $G'$ by $+i$.
Then (as we will see below)
the determinant of $\tilde{B'}$ is actually zero.
Moreover,
there is no alternative way of choosing the signs
that ensures that all the terms contribute coherently.

The problem is that
we can still take cycles in $\mu_1 + \mu_2$,
where $\mu_1$, $\mu_2$ now denote
two matchings of $G'$,
and we can still evaluate the relative parity
as $(-1)^{l+1}$ where
$2l$ is the length of the cycle,
but we can no longer evaluate the relative weight
as $(-1)^A$,
since the cycle map not be contractible
and hence does not enclose any area.
(Contractible loops pose no problem,
even though they lack a well-defined ``inside'',
since the two regions into which
a contractible cycle divides the torus
must have areas of equal parity.)
There is probably some parity-version of Pick's theorem
that applies in this setting,
at least when the edges as well as the vertices
of the polygon are required to belong to the grid,
but I haven't yet formulated, let alone proved,
such a result.

Kasteleyn says
``it can be shown that''
the determinant of $\tilde{B'}$ counts correctly
only those matchings that contain
an even number of horizontal special edges
and an even number of vertical special edges ---
that all other matchings contribute with the wrong sign.
This is in fact the case, and in private communication
Glenn Tesler, William Jockusch, and Greg Kuperberg
have sent me arguments that substantiate Kasteleyn's claim;
however, I have never seen these details anywhere in print.
Kasteleyn himself seems to have sidestepped the problem
by switching to a different approach,
in between writing of his journal article
and his later article {\it Graph Theory and Crystal Physics\/};
in the latter, he uses orientations of graphs
instead of weightings of graphs.

Taking Kasteleyn's claim as true,
we can proceed to find $N'(m,n)$.
Let us focus for simplicity on the case in which
$m$ and $n$ are both divisible by 4.
Let $\tilde{B}_0 = \tilde{B}$. 
Let $\tilde{B}_1$ be obtained from $\tilde{B}$
by changing the signs of all the $+i$'s
associated with the special vertical edges of $G'$.
Let $\tilde{B}_2$ be obtained from $\tilde{B}$
by changing the signs of all the $+1$'s
associated with the special horizontal edges of $G'$.
Let $\tilde{B}_3$ be obtained from $\tilde{B}$
by making both sorts of sign-changes.
Then it is easy to check
that the following table applies:

\[ \begin{array}{ccccc}
      & \tilde{B}_0 & \tilde{B}_1 & \tilde{B}_2 & \tilde{B}_3 \\
\mbox{(e,e)} & + & + & + & + \\
\mbox{(o,e)} & - & - & + & + \\
\mbox{(e,o)} & - & + & - & + \\
\mbox{(o,o)} & - & + & + & -
\end{array} \]
\noindent
For instance, if a matching involves
an odd number of special horizontal edges
and an even number of special vertical edges ---
that is, if the matching is of type ``(odd,even),''
or ``(o,e)'' for short ---
then it contributes $-1$ to $\det\tilde{B}_0$,
$-1$ to $\det\tilde{B}_1$,
$+1$ to $\det\tilde{B}_2$,
and $+1$ to $\det\tilde{B}_3$.
It is evident from the table
that the linear combination
$\frac{1}{2}(-\det \tilde{B}_0 + \det \tilde{B}_1 
+ \det \tilde{B}_2 + \det \tilde{B}_3)$
counts each matching with weight $+1$,
and so is equal to $N'(m,n)$.

Introducing modified adjacency matrices $\tilde{A}_k$ ($k=0,1,2,3$)
that correspond to the modified $\tilde{B}_k$'s in the obvious way,
we are able to write
\[ N'(m,n) =
\frac{1}{2}(-\sqrt{\det\tilde{A}_0} + \sqrt{\det\tilde{A}_1} 
+ \sqrt{\det\tilde{A}_2} + \sqrt{\det\tilde{A}_3}) . \]

To evaluate these determinants, introduce the matrices
$D_m^+$, $D_m^-$, $D_n^+$, and $D_n^-$,
where $D_m^+$ is the $m$-by-$m$ circulant matrix
\[ \left( \begin{array}{ccccccc}
0 & 1 & 0 & 0 & \dots & 0 & 1 \\
1 & 0 & 1 & 0 & \dots & 0 & 0 \\
0 & 1 & 0 & 1 & \dots & 0 & 0 \\
\vdots & \vdots & \vdots & \vdots & & \vdots & \vdots \\
0 & 0 & 0 & 0 & \dots & 0 & 1 \\
1 & 0 & 0 & 0 & \dots & 1 & 0 \end{array} \right) \ , \]
$D_m^-$ is the $m$-by-$m$ not-quite-circulant matrix
\[ \left( \begin{array}{ccccccc}
0 & 1 & 0 & 0 & \dots & 0 & -1 \\
1 & 0 & 1 & 0 & \dots & 0 & 0 \\
0 & 1 & 0 & 1 & \dots & 0 & 0 \\
\vdots & \vdots & \vdots & \vdots & & \vdots & \vdots \\
0 & 0 & 0 & 0 & \dots & 0 & 1 \\
-1 & 0 & 0 & 0 & \dots & 1 & 0 \end{array} \right) \ , \]
and $D_n^+$ and $D_n^-$ are similarly defined $n$-by-$n$ matrices.
We can write
$A_0 = D_m^+ \otimes I_n + I_m \otimes iD_n^+$,
$A_1 = D_m^+ \otimes I_n + I_m \otimes iD_n^-$,
$A_2 = D_m^- \otimes I_n + I_m \otimes iD_n^+$, and
$A_3 = D_m^- \otimes I_n + I_m \otimes iD_n^-$,
so we can find the eigenvalues of all four matrices
provided we can find the eigenvalues of the $D^+$ and $D^-$ matrices.
This is easily done:
the circulant matrix $D_m^+$
has eigenfunctions $x \mapsto e^{2j\pi i x/m}$
with eigenvalues $2 \cos \frac{2j \pi}{m}$
($1 \leq j \leq m$),
while the near-circulant matrix $D_m^-$
has eigenfunctions $x \mapsto e^{(2j-1) \pi i x/m}$
with eigenvalues $2 \cos \frac{(2j-1) \pi}{m}$
($1 \leq j \leq m$).
Using these values,
it can be checked that $\det A_0$ actually vanishes,
so we need not worry about potential difficulties
arising from the fact that our linear combinination of determinants
has both plus and minus signs,
allowing for massive cancellations
that might swamp the final answer.
(It might be interesting to have a combinatorial explanation
for the vanishing of $\det B_0$,
via some sort of pairing of terms.)

When all the work is done,
is turns out that $N'(m,n)$ grows
at asymptotically the same rate
as $N(m,n)$.

 
What we really are after is the {\it entropy} of the
dimer model.
This is defined as
\[
\lim_{m,n \rightarrow \infty} \frac{1}{mn} \log N^*(m,n)
\]
where $N^* (m,n)$ is the number of different possible
$m$-by-$n$ excerpts of dimer configurations on the entire plane.
Equivalently,
we may imagine laying down dominoes on a checkerboard
such that dominoes are now allowed to straddle the boundary of the board.
We call such an arrangement an ``overtiling'' of the board.
It is not hard to show that
every overtiling of a rectangle
extends to a tiling of the plane,
so $N^*(m,n)$ is simply the number of overtilings 
of an $m$-by-$n$ rectangle.
We have $N^*(m,n) \geq N'(m,n) \geq N(m,n)$,
since the ``straight'' boundary of the rectangle 
is a special case of doubly periodic boundary conditions,
and doubly periodic boundary conditions
in turn form a special case
of arbitrary boundary conditions.

Let us digress briefly to consider
why the limit
\[
\lim_{m,n \rightarrow \infty} \frac{1}{mn} \log N^*(m,n)
\]
exists.
This is a fairly straightforward generalization
of the one-dimensional argument,
but it's worthwhile checking
that the same analysis goes through.
(Thanks to Boris Solomyak 
for helping me work this out.)

We need to note that the function $N^*(\cdot,\cdot)$
is monotone in each of its arguments;
that is, for all $m' \geq m$ and $n' \geq n$,
$N^*(m',n)$ and $N^*(m,n')$ are both $\geq N(m,n)$.
We also need to note that $N^*(m,n)$
is submultiplicative in each of its arguments;
that is, if $m=m_1+m_2$
and $n=n_1+n_2$
then $N^*(m,n) \leq N^*(m_1,n) N^*(m_2,n)$
and $N^*(m,n) \leq N^*(m,n_1) N^*(m,n_2)$.
Equivalently,
$\log N^*(m,n)$ is subadditive in each argument,
so that for instance
$\log N^*(jm,kn) \leq jk \log N^*(m,n)$.

Now let
\[
\alpha = \sup_{m_0,n_0} \inf_{m \geq m_0; n \geq n_0} 
\frac{1}{mn} \log N^*(m,n) \ .
\]
Observe that 
$\frac{1}{m n} \log N^*(m_0,n_0) \geq \alpha$
for all $m,n$;
for if it happened that
$\frac{1}{m n} \log N^*(m,n) < \alpha$
then we would necessarily have
$\frac{1}{m' n'} \log N^*(m',n') < \alpha$
for arbitrarily large values of $m'$ and $n'$
(namely, those that are multiples of $m$ and $n$, respectively),
contradicting our definition of $\alpha$.

Fix $\epsilon > 0$,
and take $m_0, n_0$ such that
$\alpha \leq \frac{1}{m_0 n_0} \log N^*(m_0,n_0) < \alpha + \epsilon$.
Take $m$, $n$ suitably large
(just how large they need to be will be determined shortly)
and write $m=jm_0+r$ and $n=kn_0+s$
with $0 \leq r < m_0$, $0 \leq s < n_0$.
Since $N^*(\cdot,\cdot)$ is monotone,
$\log N^*(m,n)$ is less than or equal to $\log N^*((j+1)m_0,(k+1)n_0)$,
which by subadditivity is at most
$(j+1)(k+1) \log N^*(m_0,n_0)$.
Hence
\begin{eqnarray*}
\frac{1}{m n} \log N^*(m,n) 
& \leq & \frac{j+1}{m} \frac{k+1}{n} \log N^*(m_0, n_0) \\
& \leq & \frac{j+1}{jm_0} \frac{k+1}{kn_0} m_0 n_0 (\alpha + \epsilon) \\
& = & \frac{j+1}{j} \frac{k+1}{k} (\alpha + \epsilon) .
\end{eqnarray*}
By taking $m,n$ large,
we force $j,k$ to be large,
which 
forces the preceding expression to be
less than $\alpha + 2\epsilon$, say.
Since $\epsilon$ was arbitrary,
we have shown the existence of the limit.

We now face the truly interesting question:
What is $\alpha$?
Our earlier work tells us that
$\alpha \geq G/\pi$,
but it does not give us an upper bound.
We may write
$N^*(m,n) = \sum_C N_C^*(m,n)$,
where $C$ stands for some boundary-configuration
(specifically, the locations of those dominoes
that straddle the boundary of the rectangle)
and $N_C^*(m,n)$ is the number of ways of
tiling the rectangle
subject to the boundary condition $C$.
We know that $N^*(m,n)$ grows quadratic-exponentially,
namely as $\exp \alpha mn$,
while the number $B$ of boundary-conditions $C$
grows only linear-exponentially
(in fact, it is bounded above by $2^{4n}$).
Hence the {\it average} value of $N_C^*(m,n)$,
as $C$ ranges uniformly over all boundary conditions,
also grows like $\exp \alpha mn$,
and in particular,
there must exists a two-parameter family 
of boundary conditions $C(m,n)$
such that and $N_{C(m,n)}^*(m,n)$ grows like $\exp \alpha mn$.

Unfortunately, we do not know
what the suitable conditions $C(m,n)$ to choose.
Fortunately, we do not need to,
on account of a lovely trick discovered by Greg Kuperberg.
Picture the $m$-by-$n$ rectangle $[0,m] \times [1,n]$ in $\R^2$,
with the (unknown) boundary condition $C(m,n)$.
By reflecting $C(m,n)$ in the lines $x=m$ and $y=n$,
we get doubly-periodic boundary conditions on
$[0,2m] \times [1,2n]$,
which must be compatible with exactly
$(N_{C(m,n)}^* (m,n))^4$
overtilings of the $2m$-by-$2n$ rectangle.
Hence $N^*(2m,2n) \geq (N_{C(m,n)}^* (m,n))^4$,
Hence $\frac{1}{m} \frac{1}{n} \log (N_{C(m,n)}^* (m,n))
\leq \frac{1}{2m} \frac{1}{2n} \log N^*(2m,2n)$. 
Sending $m,n$ to infinity,
we get $\alpha \leq G/\pi$.
Combining this with the reverse inequality,
we find that $\alpha = G/\pi$.

\end{document}